\documentclass[12pt, twoside, leqno]{article}
\usepackage{amsmath,amsthm}
\usepackage{amssymb,latexsym}
\usepackage{enumerate}
\usepackage[latin1]{inputenc}

\newtheorem{defin}{Definition}[section]
\newtheorem{teore}[defin]{Theorem}
\newtheorem{lem}[defin]{Lemma}
\newtheorem{corol}[defin]{Corollary}
\newtheorem{osser}[defin]{Remark}
\newtheorem{propo}[defin]{Proposition}

\newenvironment{theorem}{\begin{teore}\sl}{\end{teore}}

\newenvironment{remark}{\begin{osser}\rm}{\end{osser}}

\newenvironment{lemma}{\begin{lem}\sl}{\end{lem}}

\newenvironment{corollary}{\begin{corol}\sl}{\end{corol}}

\newenvironment{proposition}{\begin{propo}\sl}{\end{propo}}

\newenvironment{dimos}
{\paragraph{Proof}}
{\qed}

\numberwithin{equation}{section}

\newcommand\uT{{\underline T}}
\newcommand\uX{{\underline X}}
\newcommand\uY{{\underline Y}}

\newcommand\ual{{\underline\alpha}}
\newcommand\ualpha{{\underline{\!\!\phantom{\gamma}\alpha}}}
\newcommand\uAlpha{{\underline\Alpha}}
\newcommand\ubeta{{\underline\beta}}
\newcommand\uBeta{{\underline\Beta}}
\newcommand\ugamma{{\underline\gamma}}
\newcommand\uGamma{{\underline\Gamma}}
\newcommand\utheta{{\underline{\!\!\phantom{\gamma}\theta}}}
\newcommand\uTheta{{\underline\Theta}}

\newcommand\PP{{\mathbb P}}
\newcommand\Q{{\mathbb Q}}
\newcommand\C{{\mathbb C}}
\newcommand\K{{\mathbb K}}
\renewcommand\L{{\mathbb L}}
\newcommand\R{{\mathbb R}}

\newcommand\A{{\mathbb A}}
\newcommand\Z{{\mathbb Z}}

\newcommand\DD{{\mathcal D}}
\newcommand\LL{{\mathcal L}}
\newcommand\norm{{\mathcal N}}

\newcommand\CC{{\mathcal C}}

\newcommand\ord{{\mathrm{ord}}}

\newcommand\genus{{\mathbf g}}

\newcommand\tilg{{\widetilde g}}
\newcommand\tily{{\widetilde y}}
\newcommand\tilY{{\widetilde Y}}

\newcommand\Alpha{{\mathrm A}}
\newcommand\Beta{{\mathrm B}}
\newcommand\Mu{{\mathrm M}}
\newcommand\ph\varphi

\newcommand\height{{\mathrm h}}
\newcommand\kps{{\mathrm{KPS}}}

\newcommand\hatualpha{{\widehat{\underline{\!\!\phantom{\gamma}\alpha}}}}
\newcommand\hatbeta{{\widehat{\beta}}}
\newcommand\hatubeta{{\widehat{\underline\beta}}}
\newcommand\hatgamma{{\widehat{\gamma}}}
\newcommand\hatugamma{{\widehat{\underline\gamma}}}
\newcommand\hatdelta{{\widehat{\delta}}}

\newcommand\hatutheta{{\widehat{\underline{\!\!\phantom{\gamma}\theta}}}}
\newcommand\hatphi{{\widehat{\ph}}}

\newcommand\hatCC{{\widehat \CC}}
\newcommand\hatd{{\widehat d}}
\newcommand\hatf{{\hat{f}}}
\newcommand\hatP{{\widehat P}}
\newcommand\hatx{{\widehat x}}
\newcommand\haty{{\widehat{y}}}

\newcommand\tilhaty{{\widetilde\haty}}

\begin{document}




\title{Quantitative Riemann Existence\\ Theorem over a Number Field}

\author{Yuri F.~Bilu\\
Universit\'e Bordeaux 1, Institut de Math\'ematiques\\
351 cours de la Libération, 33405 Talence, France\\
E-mail: yuri@math.u-bordeaux1.fr\and
Marco Strambi\\
Dipartimento di Matematica,
Universit\`a di Pisa\\
Lago Bruno Pontecorvo 5, 56127 Pisa, Italie\\
E-mail: strambi@mail.dm.unipi.it} 

\setcounter{tocdepth}1

\date{July 29, 2009}

\maketitle


\renewcommand{\thefootnote}{}



\renewcommand{\thefootnote}{\arabic{footnote}}
\setcounter{footnote}{0}


\hfuzz 3pt

\maketitle


\begin{abstract}
Given a covering of the projective line with ramifications over~$\bar\Q$, we define a plain model of the algebraic curve realizing the Riemann existence theorem for this covering, and bound explicitly the defining equation of this curve and its definition field. 
\end{abstract}

\section{Introduction}

The \textsl{Riemann Existence Theorem} asserts that every compact Riemann surface is (analytically isomorphic to) a complex algebraic curve. In other words, if~$f$ is a non-constant meromorphic function on a compact Riemann surface~$S$, then the field of all meromorphic functions on~$S$ is a finite extension of $\C(f)$. 

One of the most common ways of defining Riemann surfaces is realizing them as finite ramified coverings of the Riemann sphere $\PP^1(\C)$. Moreover, even if the covering is purely topological, the $\C$-analytic structure on the Riemann sphere lifts, in a unique way,  to the covering surface. Thus, the Riemann Existence Theorem can be restated as follows.

\begin{theorem}
\label{tter}
Let~$\Mu$ be a finite subset of $\PP^1(\C)$. Then for any finite covering of ${\PP^1(\C)}$ by a closed oriented surface, unramified outside the set~$\Mu$, there exists a complex algebraic curve~~$\CC$ and a rational function ${x\in \C(\CC)}$ such that our covering is isomorphic\footnote{Two morphisms of ${S_1\stackrel{\pi_1}\to S}$ and ${S_2\stackrel{\pi_2}\to S}$ of topological spaces are \textsl{isomorphic} if there exists a homeomorphism ${S_1\stackrel{\ph}\to S_2}$ such that ${\pi_1=\pi_2\circ\ph}$.} to   ${\CC(\C)\stackrel{x}\to \PP^1(\C)}$, the covering defined by~$x$. Moreover, the couple $(\CC,x)$ is unique up to a naturally defined isomorphism\footnote{If $(\CC',x')$ is another such couple, then the field isomorphism ${\C(x)\to \C(x')}$ given by ${x\mapsto x'}$, extends to a field isomorphism ${\C(\CC)\to \C(\CC')}$.}.  
\end{theorem}

We refer to~\cite{Deb01} for several more precise statements. 

The purpose of this article is to give an effective description of the curve~$\CC$, or, more precisely, of the couple $(\CC,x)$, in terms of the degree of the initial topological covering and the set~$\Mu$ of the ramification points, provided the points from that set are defined over the field~$\bar\Q$ of all algebraic numbers. In this case the curve~$\CC$ is also defined over~$\bar\Q$ (this is the ``easy'' direction of the Theorem of Belyi). We produce a plane model of~$\CC$ over~$\bar\Q$, such that one of the coordinates is~$x$, and we give explicit bounds for the degree and the height of the defining equation of this model, and of the degree and discriminant of the number field over which this model is defined. 

Notice that we do not produce a new proof of the Riemann Existence Theorem.  In fact, we do use both the existence and the uniqueness statements of Theorem~\ref{tter}. 

Let us state our principal result.  By the \textsl{height} everywhere in this article we mean the \textsl{logarithmic affine height}; see Subsection~\ref{sbohe}.

\begin{theorem}
\label{tmain}
Let ${S\to\PP^1(\C)}$ be a finite covering of degree ${n\ge 2}$ by a closed oriented surface~$S$ of genus~$\genus$, unramified outside a finite set ${\Mu\subset \PP^1(\bar\Q)}$. Put\footnote{A pedantic reader may complain that the definition of~$h$ below is formally incorrect, because  $\height(\cdot)$ is the \textsl{affine} height, and~$\Mu$ is a subset of the \textsl{projective} line. Of course, this  can be easily overcome, for instance by writing ${\PP^1=\A^1\cup\{\infty\}}$ and defining ${\height(\infty)=0}$.}
$$
\K=\Q(\Mu), \qquad h=\max\{\height(\alpha): \alpha \in \Mu\}, \qquad \Lambda=\bigl(2(\genus+1)n^2\bigr)^{10\genus n+12n}.
$$
Then there exist a number field~$\L$, containing~$\K$, an algebraic curve~~$\CC$ defined over~$\L$ and  rational functions ${x,y\in \L(\CC)}$ such that ${\L(\CC) =\L(x,y)}$ and the following is true. 

\begin{enumerate}
\item
The covering ${\CC(\C)\stackrel{x}\to \PP^1(\C)}$, defined by~$x$, is isomorphic to the given covering ${S\to\PP^1(\C)}$. 

\item 
The rational functions ${x,y\in \L(\CC)}$ satisfy the equation ${f(x,y)=0}$, where ${f(X,Y)\in \L[X,Y]}$ is an absolutely irreducible polynomial and
\begin{equation}
\label{epoly}
\deg_Xf =\genus+1, \qquad \deg_Yf=n, \qquad \height(f) \le \Lambda (h+1).
\end{equation}

\item The degree and the  discriminant of~$\L$ over~$\K$ satisfy
\begin{equation}
\label{efield}
[\L:\K]\le \Lambda, \qquad \frac{\log\norm_{\K/\Q}\DD_{\L/\K}}{[\L:\Q]}\le \Lambda (h+1),
\end{equation}
where $\norm_{\K/\Q}$ is the norm map. 

\end{enumerate}

\end{theorem}

The principal motivation of this theorem lies in the field of effective Diophantine analysis, where the covering technique is widely used. It happens quite often that only the degree of the covering and the ramification points are known, and to work with the covering curve, one needs to have an effective description of it. In particular, in~\cite{BSS10} we use Theorem~\ref{tmain} to get a user-friendly version of the Chevalley-Weil theorem, one of the main tools of Diophantine analysis.

In brief, our method of proof is as follows. First, we use the existence part of Theorem~\ref{tter} to show the existence of~$\CC$ and~$x$. Next, we define ``quasi-canonically'' a generator~$y$ of $\bar\Q(\CC)$ over $\bar\Q(x)$, and denote by $f(X,Y)$ the irreducible polynomial satisfying ${f(x,y)=0}$. Further, we show that the coefficients of this polynomial satisfy certain system of algebraic equations and inequalities, and we use the uniqueness part of Theorem~\ref{tter} to show that the system has finitely many solutions. (To be more precise, the coefficients of~$f$ form only a part of the variables involved in the equations and inequalities.)  Using this, we estimate the height of the polynomial, and the degree and discriminant of number field generated by its coefficients. 

This argument is inspired by the work of Zverovich~\cite{Zv87}, who applies rather similar approach, though he works only in the complex domain. The system of equation considered by Zverovich is simpler than ours, but we could not understand one key point in his proof of the finiteness of the number of solutions. See more on this in Section~\ref{szv}.

Our result is sensitive only to the set~$\Mu$ of ramification points, and the degree~$n$ of the covering. It would be interesting to obtain a more precise result, which depends on the more subtle elements of the ``covering data'', like the monodromy permutations associated to every ramification point. Probably, the ``correct'' statement of Theorem~\ref{tmain} must involve the notion of the Hurwitz space associated to the given topological covering, see~\cite{De01}.  Another interesting problem is to characterize our curve not in terms of the defining equation, but in more invariant terms, for instance, to estimate its Faltings height. 

In our result, the quantity~$\Lambda$ depends exponentially on~$n$. This improves on Theorem~3A from~\cite{Bi93}, where the dependence is double exponential. There are strong reasons to believe that the ``correct'' estimate is polynomial in~$n$. Indeed, this is case for a similar problem over a function field, see the recent work of Edixhoven et al.~\cite{EJS08}.

In Sections~\ref{sbohe},~\ref{spose} and~\ref{slem} we collect various auxiliary facts needed for the proof of Theorem~\ref{tmain}. The proof itself occupies Sections~\ref{seret}--\ref{sram}. In Section~\ref{szv} we very briefly discuss the work of Zverovich.

\subsection{Notation and Conventions}
\label{ssnoca}
If $F(X)$ is a polynomial in~$X$ over some field (or integral domain), and~$\beta$ is an element of this field (or domain), then we denote by $\ord_{X=\beta}F$ the order of vanishing of~$F$ at~$\beta$. Sometimes we write simply $\ord_\beta$ or even~$\ord$, when this does not lead to a confusion. We employ the same notation not only to polynomials, but also to formal  power series in ${X-\beta}$. 

We denote by~$\alpha$ the finite point $(\alpha:1)$ of the projective line~$\PP^1$, and by~$\infty$ the infinite point $(1:0)$. 

More specific notation will be introduced at the appropriate places.

\section{Heights and Algebraic Equations}
\label{sbohe}

Let ${\ual =\left( \alpha_1,\ldots,\alpha_N\right) \in \bar{\Q}^N}$ be a point with algebraic coordinates 
in the affine space of dimension~$N$.
Let~$\K$ be a number field containing $\alpha_1,\ldots,\alpha_N$ and 
$M_\K$ the set of its valuations.
We assume that every valuation $v\in M_{\mathbb{K}}$  is normalized so
that its restriction to $\Q$ is the standard infinite or $p$-adic 
valuation.
Also,  we let $\K_{v}$ be the $v$-adic completion
of $\K$, (then, in the case of an infinite $v$, the field $\K_v$ is either $\R$ or $\C$). For ${v\in M_\K}$ we put 
$$
|\ual|_v = \max\left\{|\alpha_1|_v, \ldots, |\alpha_N|_v\right\}
$$
We now define the \textsl{absolute logarithmic affine height} (in the sequel simply \textsl{height}) of the point~$\ual$ as
\begin{equation}
\label{eheight}
\height(\ual)=\frac{1}{[\K:\Q]}\sum_{v \in M_{\K}} \left[\K_{v}:\Q_v\right]
\log^+ |\ual|_v,
\end{equation}
where ${\log^+x:=\log\max\{1,x\}}$. 
It is well-known and easy to verify that the expression on the right is independent of the choice of the field~$\K$. 
The height of ${\alpha \in \bar{\Q}}$ is,  by definition, the height of the one-dimensional vector $(\alpha)$.

For a polynomial~$f$ with algebraic coefficients we denote by  $\height(f)$ the height of the vector of its coefficients, ordered somehow. More generally, the height ${\height(f_1, \ldots, f_s)}$ of a finite system of polynomials is, by definition, the height of the vector formed of all the non-zero coefficients of all these polynomials.

\subsection{Estimates for Sums and Products of Polynomials}

The following is an immediate consequence of Lemma~1.2 from~\cite{KPS01}.

\begin{lemma}
\label{lsombrapol}
Let ${f_1,\ldots,f_s}$ be polynomials in ${\bar\Q[X_1,\ldots,X_n]}$ and put 
$$
\textstyle
d=\max\{\deg f_1, \ldots, \deg f_s\}, \qquad 
h=\height(f_1, \ldots,f_s).
$$ 
Let also~$g$ be a polynomial in $\bar\Q[Y_1,\ldots,Y_s]$. Then
\begin{enumerate}
\item 
\label{iprod}
$\height\left(\prod_{i=1}^s f_i\right)\le \sum_{i=1}^s \height\left(f_i\right) + \log (n+1)\sum_{i=1}^{s-1}\deg f_i$,
\item 
\label{ig}$\height\bigl(g\left(f_1,\ldots,f_s\right)\bigr) \le \height(g) + \bigl(h+\log(s+1)+d\log(n+1)\bigr) \deg g$. \qed
\end{enumerate}
\end{lemma}

\begin{remark}
\label{rsost}
Item~(\ref{ig}) of Lemma~\ref{lsombrapol} extends to a slightly more general situation, when the polynomial~$g$ depends, besides ${Y_1, \ldots, Y_s}$ on some other indeterminates ${T_1, \ldots, T_r}$, but one substitute new polynomials only instead of the $Y_i$-s, leaving the $T_j$-s intact. In this case we again have the estimate
$$
\height\bigl(g\left(f_1,\ldots,f_s, T_1, \ldots, T_r\right)\bigr) \le \height(g) + \bigl(h+\log(s+1)+d\log(n+1)\bigr) \deg_\uY g
$$
(independently of~$r$ and of $\deg_\uT g$). Indeed, we can write ${g=\sum_k g_k(\uY)h_k(\uT)}$, where $h_k(\uT)$ are pairwise dinstct monomials in ${\uT=(T_1, \ldots, T_r)}$, and apply Lemma~\ref{lsombrapol}~(\ref{ig}) to each~$g_k$. 
\end{remark}

Here is a particular case of  Lemma~\ref{lsombrapol}, where a slightly sharper estimate holds (see~\cite{KPS01}, end of Subsection~1.1.1). 

\begin{lemma}
\label{lsombramat}
Let ${\left(f_{ij}\right)_{ij}}$ be an ${s\times s}$ matrix of polynomials in $\bar\Q[X_1,\ldots,X_n]$ of degrees and heights bounded by~$d$ and~$h$, respectively. Then
$$
\height\bigl(\det\left(f_{ij}\right)_{ij}\bigr) \le s\bigl(h+\log s+d\log(n+1)\bigr). \eqno\square
$$
\end{lemma}

We need one more technical lemma. 

\begin{lemma}
\label{lrho}
Let ${g(X,Y)\in \bar\Q[X,Y]}$ be of $X$-degree~$m$, and fix ${\rho\in \bar\Q}$. Put 
$$
f(X,Y):=(X-\rho)^m g\bigl((X-\rho)^{-1},Y\bigr).
$$
Then
$$
\height(f)\le \height(g)+m\height(\rho)+ 2m\log 2.
$$
\end{lemma}

\begin{dimos}
The polynomials $g(X,Y)$ and ${\tilg(X,Y):=X^mg(X^{-1},Y)}$ have the same coefficients and thereby the same height. Applying Lemma~\ref{lsombrapol} and Remark~\ref{rsost}, we obtain the result. 
\end{dimos}

\subsection{Bounds for Solutions of Algebraic Equations}
By an \textsl{algebraic set}  we mean a subset of $\bar\Q^N$, defined by a system of polynomial equations. We treat algebraic sets as in \cite[16. Kapitel]{vW93} (where they are called \textit{algebraische Mannigfaltigkeiten}), that is,   purely set-theoretically, without counting multiplicities. By a \textsl{component} of an algebraic sets we mean an irreducible component.

Let ${p_1(\uX), \ldots, p_k(\uX)}$ be polynomials in ${\uX=(X_1, \ldots, X_N)}$ 
with algebraic coefficients. By an   \textsl{isolated 
solution} of the system of polynomial equations 
\begin{equation}
\label{esys}
p_1(\uX)=\ldots =p_k(\uX)=0.
\end{equation}
we mean  a zero-dimensional component   of the algebraic set in $\bar\Q^N$ defined by~(\ref{esys}). (Existence of such a component implies that ${k\ge N}$.) Our aim is to 
bound the height of an isolated solution in terms of the degrees and the heights of the polynomials ${p_1, \ldots, p_k}$. 

Such a bound follows from the arithmetical Bézout inequality due to Bost, Gillet and Soulé~\cite{BGS94} and Philippon~\cite{Ph91}. Krick, Pardo and Sombra~\cite{KPS01} did a great job of producing a user-friendly version of this fundamental result. We very briefly recall some facts from~\cite{KPS01} which will be used here. 
For an affine algebraic set ${V\subset \A^N}$, defined over~$\bar\Q$, Krick, Pardo and Sombra \cite[Section~1.2]{KPS01} define the \textsl{height} of~$V$, to be denoted here as  $\height_\kps(V)$. We do not reproduce here the full definition of this height function, but only list four of its properties. 

\begin{proposition}
The Krick-Pardo-Sombra height function has the following properties.

\begin{description}

\item[(positivity)] For any~$V$ we have ${\height_\kps(V)\ge 0}$.

\item[(additivity)]
The height function is ``additive'' in the following sense: for any~$V_1$ and~$V_2$ without common components we have 
$$
\height_\kps(V_1\cup V_2)=\height_\kps(V_1)+ \height_\kps(V_2).
$$

\item[(one-point set)]
If $\{\ual\}$ is a one-point algebraic set, then 
${\height(\ual)\le \height_\kps(\{\ual\})}$

\item[(Bézout inequality)]
Let~$V$ be the algebraic set defined by 
$$
p_1(\uX)=\cdots=p_N(\uX)=0,
$$
where  ${p_i(\uX) \in \bar\Q(\uX)}$
for ${i=1\ldots,N}$.
Put
\begin{equation}
\label{ensh}
\begin{gathered}
\nabla = \deg p_1\cdots \deg p_N, \qquad \Sigma=\frac{1}{\deg p_1}+\cdots+\frac{1}{\deg p_N}, \\
h=\max\{\height(p_1), \ldots, \height(p_N)\}.
\end{gathered}
\end{equation}
Then
\begin{equation}
\label{ehkps}
\height_\kps(V) \le \nabla\Sigma h +2 \nabla N\log (N+1).
\end{equation}

\end{description}
\end{proposition}

\begin{dimos}
Positivity and additivity follow immediately from the definition. For the height of a one-point set see \cite[end of Section~1.2.3]{KPS01}; in fact, $\height_\kps(\{\ual\})$ is defined as the right-hand side~(\ref{eheight}) but with ${\log^+ |\ual|_v}$ replaced by ${\log\left(1+ |\alpha_1|_v^2+ \cdots+ |\alpha_N|_v^2\right)^{1/2}}$ for archimedean~$v$.
Finally, for the Bézout inequality see Corollary~2.11 from~\cite{KPS01}, or, more precisely,   the displayed inequality just before the beginning of Section~2.2.3 on page 555 of~\cite{KPS01}.
\end{dimos}

\medskip
We adapt the work of Krick, Pardo and Sombra  as follows.

\begin{proposition}
\label{pbound}
Let~$K$ be a number field and let ${p_1(\uX), \ldots, p_k(\uX)\in \K[\uX]}$ be polynomials in ${\uX=(X_1, \ldots, X_N)}$. Let~$\ual$ be an isolated solution of~(\ref{esys}) and ${\L=\K(\ual)}$ the number field generated by the coordinates of~$\ual$. Then ${k\ge N}$. Further, assume that
$$
\deg p_1\ge \deg p_2\ge \ldots\ge \deg p_k.
$$
Let also~$\nabla$, $\Sigma$ be defined as in~(\ref{ensh}) and and ${h=\max\{\height(p_1), \ldots, \height(p_k)\}}$.
Then 
\begin{align}
\label{edegl}
[\L:\K]&\le \nabla,\\ 
\label{ehea}
[\L:\K]\height(\ual)&\le \nabla\Sigma h +2 \nabla N\log (N+1)  ,\\
\label{edisl}
\frac{\log\norm_{\K/\Q}\DD_{\L/\K}}{[\L:\Q]}&\le 2\nabla\Sigma h +5 \nabla N\log (N+1) ,
\end{align}
where $\DD_{\L/\K}$ is the discriminant of~$\L$ over~$\K$ and $\norm_{\K/\Q}$ is the norm map. 
\end{proposition}

The following consequence is immediate. 

\begin{corollary}
\label{cbound}
In the set-up of Proposition~\ref{pbound}, denote by~$V$ the algebraic subset of~$\bar\Q^N$ defined by~(\ref{esys}), and let~$W$ be another algebraic subset of~$\bar\Q^N$ such that the difference set ${V\setminus W}$ is finite. 
Then every ${\ual\in V\setminus W}$ satisfies~(\ref{edegl}),~(\ref{ehea}) and~(\ref{edisl}). \qed
\end{corollary}

For the proof of Proposition~\ref{pbound} we shall use the following lemma, due to Silverman \cite[Theorem~2]{Si84}.
\begin{lemma}
\label{lsil}
Let~$\K$ be a number field and~$\ual$ be a point in~$\bar\Q^N$. 
Then the relative discriminant $\DD_{\L/\K}$ of the field ${\L=\K(\ual)}$ over~$\K$ satisfies
$$
\frac{\log\norm_{\K/\Q}\DD_{\L/\K}}{[\L:\Q]}\le 2([\L:\K]-1)\height(\ual) +\log [\L:\K]. \eqno\square
$$
\end{lemma}

\paragraph{Proof of Proposition~\ref{pbound}}
We denote by~$V$ the algebraic set defined by~(\ref{esys}). Since it has a $0$-dimensional component~$\ual$, we have ${k\ge N}$. 
Among the~$k$ polynomials ${p_1, \ldots, p_k}$ one can select~$N$ polynomials ${q_1, \ldots, q_N}$ such that~$\ual$ is an isolated solution of the system ${q_1(\uX)=\ldots=q_N(\uX)}$. The algebraic set defined by this system has at most ${\deg q_1\cdots \deg q_N\le \nabla}$ irreducible (over~$\bar\Q$) components: this follows from the geometric Bézout inequality. In particular,  there is at most~$\nabla$ isolated solutions.
Since a $\K$-conjugate of an isolated solution  is again an isolated solution, we must have~(\ref{edegl}). 

Further, the four properties above imply that 
\begin{equation}
\label{ehadd}
\sum_{\text{$\{\ual\}$ component of~$V$}}\height(\ual) \le 
\height_\kps(V)\le\nabla\Sigma h +2 \nabla N\log (N+1),              
\end{equation}
where the sum is over  the $0$-dimensional components of $V(\bar\Q)$. 
Since all conjugates of~$\alpha$ have the same height, the left-hand side of~(\ref{ehadd}) exceeds ${[\L:\K]\height(\ual)}$, which proves~(\ref{ehea}). Combining it with 
Lemma~\ref{lsil}, we obtain~(\ref{edisl}). \qed

\section{Power Series}
\label{spose}
In this section~$K$ is a field of characteristic~$0$ and ${f(X,Y)\in K[[X]][Y]}$ is a polynomial in~$Y$ with coefficients in the ring $K[[X]]$ of formal power series. We denote by $\ord$ the order of vanishing  at~$0$. By the \textsl{initial segment of length~$\kappa$} (or simply \textsl{$\kappa$-initial segment}) of a power series ${y=\sum_{k=0}^\infty \gamma_kX^k}$ we mean ${y=\sum_{k=0}^\kappa \gamma_kX^k}$. 

\begin{lemma}
\label{lhens}
Let ${\tily = \sum_{k=0}^\kappa \gamma_kX^k\in K[X]}$ be a polynomial in~$X$ of degree at most~$\kappa$. Assume that
$$
\qquad \ord f(X,\tily)>2\kappa, \qquad \ord f'_Y(X,\tily)=\kappa.
$$
Then  there exists a unique formal power series ${y=\sum_{k=0}^\infty \gamma_kX^k\in K[[X]]}$ such that ${f(X,y)=0}$, and such that~$\tily$ is the initial segment of~$y$ of length~$\kappa$.
\end{lemma}

\begin{dimos}
By Hensel's Lemma, there exists a unique power series~$y$ such that ${f(X,y)=0}$ and ${\ord(y-\tily)>\kappa}$. The latter inequality implies that~$\tily$ is the initial segment of~$y$ of length~$\kappa$. 
\end{dimos}

{\sloppy

\begin{lemma}
\label{lhens1}
\begin{enumerate}

\item
\label{ihen}
Let~${y\in K[[X]]}$ be a formal power series such that 
${f(X,y)=0}$.
We define  ${\kappa=\ord f'_Y(X,y)}$ and we let~$\tily$ be the initial segment of~$y$ of length~$\kappa$.
Then
${\ord f(X,\tily) > 2\kappa}$ and  ${\ord f'_Y(X,\tily) =\kappa}$.

\item
\label{idist}
Let ${y_1,y_2\in K[[X]]}$ be distinct formal power series such that 
$$
f(X,y_1)=f(X,y_2)=0,
$$ 
and let $\kappa_1,\kappa_2$ be defined as~$\kappa$ in part~(\ref{ihen}). 
Then   the $k$-th coefficients of~$y_1$ and~$y_2$ are distinct for some  ${k\le \min \{\kappa_1, \kappa_2\}}$. 
\end{enumerate}
\end{lemma}

}

\begin{dimos}
Since~$\tily$ is the $\kappa$-initial segment of~$y$, we have ${\ord (y-\tily)>\kappa}$. Hence
$$
f(X,\tily)= f(X,y)+f'_Y(X,y)(y-\tily)+ \text{terms of order $>2\kappa$},
$$
Since ${f(X,y)=0}$ and ${\ord f'_Y(X,y)=\kappa}$, the right-hand side is of order $>2\kappa$. Similarly, 
$$
f'_Y(X,\tily)= f'_Y(X,y)+ \text{terms of order $>\kappa$},
$$
which implies that the right-hand side is of order~$\kappa$. We have proved part~(\ref{ihen}). 

Now to part~(\ref{idist}). Lemma~\ref{lhens} implies that~$y_j$ is the single power series satisfying ${f(X,y_j)=0}$ and having~$\tily_j$ as an initial segment. Since the series~$y_1$ and~$y_2$ are distinct, none of~$\tily_j$ can be an initial segment of the other\footnote{If, say,~$\tily_1$ is an initial segment of~$\tily_2$ then the same argument as above shows that ${\ord f'_Y(X,\tily_2)=\ord f'_Y(X,\tily_1)}$, that is, ${\kappa_1=\kappa_2}$, whence ${\tily_1=\tily_2}$. Lemma~\ref{lhens} now implies that ${y_1=y_2}$, a contradiction.}. Whence the result. 
\end{dimos}


\begin{lemma}
\label{lunram}
Suppose~$K$ algebrically closed and let~${y_1,\ldots, y_\ell\in K[[X]]}$ be pairwise distinct formal power series such that 
$$
f(X,y_1)=\ldots=f(X,y_\ell)=0.
$$
Assume that the polynomial~$f$ is monic in~$Y$ (that is,~$f$ is of the form ${Y^n+\text{termes of lower degree in~$Y$}}$) and that 
\begin{equation}
\label{esumk}
\sum_{j=1}^\ell\ord f'_Y(y_j)=\ord\, d(X),
\end{equation}
where $d(X)$ is the $Y$-discriminant of~$f$. Then~$f$ splits into linear factors over the ring $K[[X]]$:
$$
f(X,Y)=(Y-y_1)\cdots (Y-y_n),
$$
where ${y_1,\ldots, y_n\in K[[X]]}$. 
\end{lemma}

{\sloppy

\begin{dimos}
Since~$f$ is monic, it splits, by the Puiseux theorem, into linear factors over the ring $K[[X^{1/e}]]$ for some~$e$:
$$
f(X,Y)=(Y-y_1)\cdots (Y-y_n),
$$
where ${y_{\ell+1},\ldots, y_n\in K[[X^{1/e}]]}$. Further, 
${d(X)=\prod_{j=1}^n f'_Y(y_j)}$, 
which, together with~(\ref{esumk}) implies that 
\begin{equation}
\label{eorder}
\ord f'_Y(y_j)=0 \qquad (j=\ell+1, \ldots, n). 
\end{equation}
If we now write ${y_j= a_{j0}+a_{j1}X^{1/e}+\ldots}$, then~(\ref{eorder}) implies that
$$
\ord f'_Y(X,a_{j0})=0 \qquad (j=\ell+1, \ldots, n). 
$$
Lemma~\ref{lhens} now implies that in each of the rings $K[[X]]$ and $K[[X^{1/e}]]$, the polynomial~$f$ has exactly one root with initial term~$a_{j0}$. Hence ${y_j\in K[[X]]}$ for  ${j=\ell+1, \ldots, n}$, as wanted. 
\end{dimos}

}

\section{Miscellaneous Lemmas} 
\label{slem}


\begin{lemma}
\label{lpoles}
Let~$\CC$ be a smooth projective curve defined over an algebraically closed field~$K$ of characteristic~$0$. Let ${x\in K(\CC)}$ have only simple poles, and let ${y\in K(\CC)}$ have a single (possibly, multiple) pole which is a pole of~$x$ as well. Then ${K(\CC)=K(x,y)}$.
\end{lemma}

\begin{dimos}
Since~$x$ has only simple poles in $K(\CC)$, the place at~$\infty$ of the field $K(x)$ splits completely in $K(\CC)$. Let~$P$ be the pole of~$y$, and let~$\widetilde{P}$ be the place of $K(x,y)$ below~$P$. Then~$\widetilde{P}$ is above the place at~$\infty$ of $K(x)$. Hence it also splits completely in~$K(\CC)$.  

Now assume that $K(x,y)$ is a proper subfield of $K(\CC)$. Then there are at least~$2$ places of $K(\CC)$ above~$\widetilde{P}$. In particular, there is a place ${P'\ne P}$ above~$\widetilde{P}$. This~$P'$ must be a pole of~$y$, a contradiction.  
\end{dimos}

\begin{lemma}
\label{lber}
Let~$K$ be an algebraically closed field of characteristic~$0$ and~$V$ is a non-empty quasiprojective variety over~$K$. Let  ${\left\{(\CC_t, D_t): t\in V\right\}}$ be an algebraic family of curves supplied with an effective divisor. Also, let~$s$ be a positive integer. 
Assume that there exists ${\tau\in V}$  such that~$\CC_\tau$ is irreducible and ${\dim\LL(D_\tau)=s}$.  Then  the set 
$$
\left\{t\in V : 
\begin{array}l
\text{either~$\CC_t$ is reducible}\\ \text{or~$\CC_t$ is irreducible and ${\dim\LL(D_t) > s}$}
\end{array}
\right \}
$$ 
is not Zariski dense in~$V$.
\end{lemma}

\begin{dimos}
This is a consequence of the theorems of Bertini and semi-continuity, see, for instance, Theorem 12.8 in~\cite[Chapter~III]{Ha93}.
\end{dimos}

\begin{lemma}
\label{lfin}
Given a positive integer~$n$ and a finite set ${\Mu\subset \C}$, there exist only finitely many extensions of the rational function field $\C(x)$ of degree~$n$, unramified outside~$\Mu$.
\end{lemma}

\begin{dimos}
This lemma (which may be viewed as an analogue of the Hermite theorem for function fields)  is an immediate consequence of the uniqueness statement of Theorem~\ref{tter}. Alternatively, it is a direct consequence of the fact that the fundamental group of a compact Riemann surface is finitely generated. 
\end{dimos}


\section{Launching the Proof of Theorem~\ref{tmain}}

\label{seret}

Let ${S\to \PP^1(\C)}$ be a covering as in the statement of Theorem~\ref{tmain}. According to Theorem~\ref{tter}, our covering is isomorphic to ${\CC(\C)\stackrel x\to\PP^1(\C)}$, where~$\CC$ is a complex algebraic curve and~$x$ is a rational function on~$\CC$. Since all ramification points of the latter covering are algebraic, the curve~$\CC$ the function~$x$ are definable over~$\bar\Q$.

{\sloppy

We are going to find a number field ${\L\supset\K}$, a rational function ${y\in \L(\CC)}$  such that ${\bar\Q(\CC)=\bar\Q(x,y)}$, and an absolutely irreducible polynomial ${f(X,Y)\in \L[X,Y]}$ such that ${f(x,y)=0}$, and such that   the degrees $\deg_Xf$, $\deg_Yf$, the height $\height(f)$, as well as the degree $[\L:\K]$ and the relative discriminant of $\L/\K$ satisfy required (in)equalities. To achieve this, we define algebraic sets~$V$ and~$W$ in a high-dimensional affine space, 
such that the  set $V\setminus W$ contains a point having the coefficients of~$f$ as part of its coordinates. We then show that the set $V\setminus W$ is finite (and hence the coefficients of~$f$) using Corollary~\ref{cbound}. As a by-product, we will also bound the degree and the discriminant of the field generated by the coefficients. 

}

We write
$$
\Mu=\{\alpha_1, \ldots, \alpha_\mu\}.
$$
For the main part of the proof we shall assume that the curve~$\CC$ is unramified over~$\infty$ (that is,~$\infty$ is not one of the points ${\alpha_1, \ldots, \alpha_\mu}$), and that~$\CC$ has no Weierstrass point above~$\infty$. In other words,   the poles of~$x$ are neither ramified nor Weierstrass.   The general case easily reduces to this one, see Section~\ref{sram}.

Now we start the detailed proof. Since it is going to be long and involved, we divide it into short logically complete steps.

\section{Function~$y$ and Polynomial $f(X,Y)$}
\label{sspoly}
Fix a pole~$P$ of~$x$. Since~$P$ is not a Weierstrass point of~$\CC$, we have 
$$
\dim \LL(mP)=2, \qquad \dim \LL((m-1)P)=1.
$$
with ${m = \genus(\CC) +1}$.

Since~$x$ is unramified above the infinity,~$x^{-1}$ can serve as a local parameter at~$P$. 
If~$y$ belongs to $\LL(mP)$, but not to $\LL((m-1)P)$, then~$y$ has the Puiseux expansion at~$P$  of the form 
${\sum_{k=-m}^\infty c_kx^{-k}}$ with ${c_{-m}\ne 0}$. Since ${\dim \LL(mP)=2}$, there exists a unique ${y\in \LL(mP)}$ with the properties 
\begin{equation}
\label{ecmco}
c_{-m}=1, \qquad
c_0=0.
\end{equation}
In the sequel, we mean by~$y$ the function satisfying these conditions.

The function~$y$ has a single pole~$P$ which is a pole of~$x$ as well. Lemma~\ref{lpoles} implies now that ${\bar\Q(\CC)=\bar\Q(x,y)}$ (here we use the assumption that~$x$ is unramified above~$\infty$). Also, since~$y$ has no poles outside the poles of~$x$, it is integral over the ring $\bar\Q[x]$. Hence, there exists a unique absolutely irreducible polynomial ${f(X,Y)\in \bar\Q[X,Y]}$, such that ${f(x,y)=0}$, monic in~$Y$ and satisfying
$$
\deg_Yf=[\bar\Q(\CC):\bar\Q(x)]=n.
$$
We also have 
$$
\deg_Xf= [\bar\Q(\CC):\bar\Q(y)]=\deg(y)_\infty=m,
$$
where ${(y)_\infty=mP}$ is the divisor of poles of~$y$. 
We write
\begin{equation}
\label{ecoeff}
f(X,Y)=  Y^n +\sum_{j=0}^{n-1} \sum_{i=0}^{m} \theta_{ij}X^iY^j.
\end{equation}

\section{Discriminant, its Roots, and Puiseux Expansions}
\label{spui}
Let $d(X)$ be the discriminant of $f(X,Y)$ with respect to $Y$. Every~$\alpha_i$ is a root of $d(X)$. Besides the $\alpha_i$-s, the polynomial $d(X)$ may have other roots; we denote them ${\beta_1, \ldots, \beta_\nu}$. Thus, we have
\begin{equation}
\label{edisc}
d(X)=\delta\prod_{i=1}^\mu\left(X-\alpha_i\right)^{\sigma_i}\prod_{i=1}^\nu\left(X-\beta_i\right)^{\tau_i},
\end{equation}
where ${\delta \in \bar\Q^\ast}$ and where~$\sigma_i$ and~$\tau_i$ are positive integers.

Now fix ${i\in \{i, \ldots, \nu\}}$. Since~$x$ is unramified over~$\beta_i$, the function~$y$ has~$n$ Puiseux expansions at~$\beta_i$ of the form
$$
y_{ij}=\sum_{k=0}^\infty\gamma_{ijk}\left(x-\beta_i\right)^k \qquad (j=1, \ldots, n).
$$
We put 
$$
\kappa_{ij}=\ord_{\beta_i} f'_Y\left(x,y_{ij}\right).
$$
Then
\begin{equation}
\label{ekeq}
\kappa_{i1}+\cdots+\kappa_{in}=\tau_i.
\end{equation}
We may assume that ${\kappa_{i1}\ge \ldots\ge\kappa_{in}}$ and we define~$\ell_i$ from the condition 
\begin{equation}
\label{eli}
\kappa_{i\ell_i}>0, \qquad \kappa_{ij}=0 \quad \text{for ${j>\ell_i}$}. 
\end{equation}
Then~(\ref{ekeq}) reads
\begin{equation}
\label{eskij}
\sum_{j=1}^{\ell_i}
\kappa_{ij}=\tau_i,
\end{equation}
which implies that 
\begin{equation}
\label{ekeq1}
\sum_{\genfrac{}{}{0pt}{}{1 \le i\le \nu}{1 \le j \le \ell_i}}(\kappa_{ij}+1) \le \sum_{\genfrac{}{}{0pt}{}{1 \le i\le \nu}{1 \le j \le \ell_i}}2\kappa_{ij}= 2(\tau_1+\cdots+\tau_\nu)\le 2\deg d(X). 
\end{equation}
This inequality will be used in Section~\ref{sindi}.

We also let~$\tily_{ij}$ be the initial segment of the series~$y_{ij}$ of  length $\kappa_{ij}$:
\begin{equation}
\label{etyij}
\tily_{ij}=\sum_{k=0}^{\kappa_{ij}}\gamma_{ijk}\left(x-\beta_i\right)^k.
\end{equation}
By Lemma~\ref{lhens1} we have 
$$
\ord_{\beta_i} f\left(x,\tily_{ij}\right)>2\kappa_{ij}, \quad  \ord_{\beta_i} f'_Y\left(x,\tily_{ij}\right)=\kappa_{ij}.
$$ 
Lemma~\ref{lhens1} also implies that, for every fixed~$i$, neither of ${\tily_{i1}, \ldots, \tily_{in}}$ is an initial segment of the other. In other words, for every distinct ${j_1,j_2\in\{1, \ldots, n\}}$ there exists a non-negative integer ${\lambda(i,j_1,j_2)\le\min\left\{\kappa_{ij_1}, \kappa_{ij_2}\right\}}$ such that 
$$
\gamma_{ij_1\lambda(i,j_1,j_2)}\ne\gamma_{ij_2\lambda(i,j_1,j_2)}.
$$

\section{Expansions at Infinity}
\label{spuinf}

We also have the Puiseux expansions of~$y$ at infinity:
\begin{equation}
\label{epuinf}
\begin{aligned}
y_{\infty j} &=\sum_{k=0}^\infty\gamma_{\infty jk}\, x^{-k} \qquad (j=2, \ldots, n),\\
y_{\infty 1} &=\sum_{k=-m}^\infty\gamma_{\infty 1k} \, x^{-k}.
\end{aligned}
\end{equation}
We define the polynomials 
$$
g(T,Y)=T^mf\left(T^{-1}, Y\right), \qquad h(T,Y)=T^{m(n+1)}f\left(T^{-1}, T^{-m}Y\right)
$$
and put ${t=x^{-1}}$,  so that the expansions~(\ref{epuinf}) can be written in powers of~$t$.
Now we define the numbers 
\begin{align*}
\kappa_{\infty j}&=\ord_{t=0}g'_Y\left(t,y_{\infty j}\right) \qquad (j=2, \ldots, n), \\
\kappa_{\infty 1}&=\ord_{t=0}h'_Y\left(t,t^my_{\infty 1}\right).
\end{align*}
We have ${h(T,T^mY)=T^{mn}g(T,Y)}$, whence 
$$
\kappa_{\infty 1}=mn+\ord_{t=0}g'_Y\left(t,y_{\infty 1}\right).
$$ 
Hence the sum ${\kappa_{\infty1}+\kappa_{\infty2}+\cdots+\kappa_{\infty n}}$ is bounded by $mn$ plus the order at ${T=0}$ of the $Y$-discriminant of $g(T,Y)$. Bounding the latter order by the degree of this discriminant, we obtain
\begin{equation}
\label{kinfineq}
\kappa_{\infty1}+\kappa_{\infty2}+\cdots+\kappa_{\infty n}\le mn+\deg d(X).
\end{equation}
Putting
\begin{equation}
\label{elinf}
\ell_\infty=n,
\end{equation}
we re-write~(\ref{kinfineq}) as  
\begin{equation}
\label{kinfineq1}
\sum_{1\le j\le \ell_\infty}(\kappa_{\infty j}+1)\le (m+1)n+\deg d(X).
\end{equation}
This will be used in Section~\ref{sindi}.

Further, for ${j=2, \ldots, n}$ let~$\tily_{\infty j}$ be the initial segment of the series~$y_{\infty j}$ of the length $\kappa_{\infty j}$, and let $\tily_{\infty 1}$ be the initial segment of the series $y_{\infty 1}$ of the length $\kappa_{\infty 1}$:
\begin{align}
\label{etyiinf}
\tily_{\infty j}&=\sum_{k=0}^{\kappa_{\infty j}}\gamma_{\infty jk} \,t^k \qquad (j=2, \ldots, n), \\
\label{ety1inf}
\tily_{\infty 1}&=\sum_{k=-m}^{\kappa_{\infty1}-m}\gamma_{\infty 1k} \,t^k.
\end{align}
Then we have 
\begin{align*}
\ord_{t=0} g\left(t,\tily_{\infty j}\right)&>2\kappa_{\infty j}, & \ord_{t=0} g'_Y\left(t,\tily_{\infty j}\right)&=\kappa_{\infty j}\qquad (j=2, \ldots, n), \\
\ord_{t=0} h\left(t,t^m\tily_{\infty 1}\right)&>2\kappa_{\infty 1}, & \ord_{t=0} h'_Y\left(t,t^m\tily_{\infty 1}\right)&=\kappa_{\infty 1}.
\end{align*}
Identities~(\ref{ecmco}) now become
$$
\gamma_{\infty1,-m}=1, \qquad \gamma_{\infty1\,0}=0.
$$
As in the finite case, for every distinct ${j_1,j_2\in\{2, \ldots, n\}}$ there exists a non-negative integer ${\lambda(\infty,j_1,j_2)\le\min\left\{\kappa_{\infty j_1}, \kappa_{\infty j_2}\right\}}$ such that 
$$
\gamma_{\infty j_1\lambda(\infty,j_1,j_2)}\ne\gamma_{ \infty j_2\lambda(\infty,j_1,j_2)}.
$$

\section{Indeterminates}
\label{sindi}
We consider the vector 
$$
\ph=\left(\utheta, \ualpha, \ubeta,\ugamma,\delta\right) \in \bar\Q^\Omega,
$$
where the dimension~$\Omega$ is defined below in~(\ref{eomegga}). Here:
\begin{itemize}
\item ${\utheta=\left(\theta_{ij}\right)_{\genfrac{}{}{0pt}{}{0\le i\le m}{0\le j\le n-1}}}$ is  the vector of  coefficients of~$f$, see~(\ref{ecoeff});

\item 
${\ualpha=(\alpha_i)_{1\le i\le \mu}}$ and
${\ubeta=\left(\beta_i\right)_{1 \le i\le \nu}}$ are the vectors of roots of the discriminant $d(X)$, and~$\delta$ is its leading coefficient, see~(\ref{edisc}); 

\item 
${\ugamma=\left(\ugamma_{\,ij}\right)_{\genfrac{}{}{0pt}{}{i\in \{1, \ldots, \nu, \infty\}}{1\le j\le \ell_i}}}$, where $\ell_i$ are defined in~(\ref{eli}) and~(\ref{elinf}), and $\ugamma_{\,ij}$ is the vector of coefficients of the initial segment~$\tily_{ij}$ of the Puiseux expansion~$y_{ij}$, see~(\ref{etyij}),~(\ref{etyiinf}) and~(\ref{ety1inf}); that is, 
${\ugamma_{\,ij}=\left(\gamma_{ijk}\right)_{0\le k\le\kappa_{ij}}}$
for ${(i,j)\ne (\infty,1)}$ and  
${\ugamma_{\,\infty 1}=\left(\gamma_{\infty 1k}\right)_{-m\le k\le\kappa_{\infty1}-m}}$;
\end{itemize}
We are only interested in the vectors~$\utheta$ and $\ualpha$, but we cannot study them separately of the other vectors defined above. 

The dimension $\Omega$ is defined by
\begin{equation}
\label{eomegga}
\Omega =(m+1)n+\mu+\nu+ \sum_{\genfrac{}{}{0pt}{}{1 \le i\le \nu}{1 \le j \le \ell_i}}(\kappa_{ij}+1)+ \sum_{1\le j\le \ell_\infty}(\kappa_{\infty j}+1)+1.
\end{equation}
We have
\begin{equation}
\label{eOmega}
\Omega\le2(m+1)n+ 4\deg d(X) +1\le 10mn+2n-8m+1,
\end{equation}
where we use~(\ref{ekeq1}),~(\ref{kinfineq1}) and the estimates 
${\mu+\nu\le\deg (d(X)) \le 2m(n-1)}$.

We shall define algebraic sets~$V$ and~$W$ in $\bar\Q^\Omega$ such that ${\ph\in V\setminus W}$ and 
${V\setminus W}$ is finite. This will allow us to use Corollary~\ref{cbound} to bound the height of~$\ph$. This would imply a bound on the height of~$\utheta$, which is the height of the polynomial~$f$.

To define our algebraic sets, we introduce the vector of indeterminates~$\Phi$ whose coordinates correspond to the coordinates of~$\ph$:
$$
\Phi=\left(\uTheta, \uAlpha, \uBeta,\uGamma,\Delta\right),
$$
where
$$
\uTheta=\left(\Theta_{ij}\right)_{\genfrac{}{}{0pt}{}{0\le i\le m}{0\le j\le n-1}}, \quad \uAlpha=(\Alpha_i)_{1\le i\le \mu},\quad
\uBeta=\left(\Beta_i\right)_{1 \le i\le \nu}, \quad \uGamma=\left(\uGamma_{\,ij}\right)_{\genfrac{}{}{0pt}{}{i\in\{1, \ldots, \nu, \infty\}}{1\le j\le \ell_i}}
$$
with
$$
\uGamma_{\,ij}=\left(\Gamma_{ijk}\right)_{0\le k\le\kappa_{ij}} \qquad \text{for $(i,j)\ne (\infty,1)$}, \qquad
\uGamma_{\,\infty 1}=\left(\Gamma_{\infty 1k}\right)_{-m\le k\le\kappa_{\infty1}-m}.
$$

\section{The Algebraic Set~$V$}
The first series of equations defining the algebraic set~$V$ is 
\begin{equation}
\label{eramV}
\Alpha_i=\alpha_i \qquad (i=1, \ldots, \mu). 
\end{equation}
To write down the rest of the equations defining~$V$ we introduce the polynomials $F(X,Y)$, $D(X)$, $G(T,Y)$ and $H(T,Y)$ with coefficients in $\Z[\uTheta]$,  which correspond to the polynomials $d(X)$, $g(T,Y)$ and $h(T,Y)$ from Section~\ref{spui}. More specifically, we put
$$
F(X,Y)= Y^n +\sum_{j=0}^{n-1} \sum_{i=0}^{m} \Theta_{ij}X^iY^j \in \Z[\uTheta][X,Y],
$$
we define $D(X)$ as the $Y$-discriminant of $F(X,Y)$ and we put
$$
G(T,Y)=T^m F\left(T^{-1}, Y\right), \qquad H(T,Y)=T^{m(n+1)} F\left(T^{-1}, T^{-m}Y\right).
$$
The second series of equations comes out from the equality
\begin{equation}
\label{edisV}
D(X)=\Delta\prod_{i=1}^\mu\left(X-\Alpha_i\right)^{\sigma_i}\prod_{i=1}^\nu\left(X-\Beta_i\right)^{\tau_i},
\end{equation}
where the quantities~$\sigma_i$ and~$\tau_i$ are defined in~(\ref{edisc}). 
In order to define the third set of equation we introduce the polynomials
\begin{align*}
\tilY_{ij}&=\sum_{k=0}^{\kappa_{ij}}\Gamma_{ijk}\left(X-\Beta_i\right)^k&(1\le i\le\nu,\quad 1\le j\le\ell_i),\\
\tilY_{\infty j}&=\sum_{k=0}^{\kappa_{\infty j}}\Gamma_{\infty jk}T^k & (2\le j\le\ell_\infty)
\end{align*} 
and the Laurent polynomial
$$
\tilY_{\infty 1}=\sum_{k=-m}^{\kappa_{\infty\!1}-m}\Gamma_{\infty 1 k}T^k.
$$
The equations come out from the relations
\begin{align}
\label{eserV}
&\begin{gathered}
\ord_{X=\Beta_i}F(X,\tilY_{i,j})> 2 \kappa_{ij},\\
\ord_{X=\Beta_i}F'_Y(X,\tilY_{i,j})\ge \kappa_{ij}
\end{gathered}
&&(1\le i\le\nu, \quad 1\le j\le \ell_i),\\[2mm]
\label{eserVinf1}
&\begin{gathered}
\ord_{T=0}G(T,\tilY_{\infty,j})> 2 \kappa_{\infty j},\\
\ord_{T=0}G'_Y(T,\tilY_{\infty,j})\ge \kappa_{\infty j}
\end{gathered}
&&(2\le j\le\ell_\infty),\\[2mm]
\label{eserVinf2}
&\begin{gathered}
\ord_{T=0}H(T,T^m\tilY_{\infty,1})> 2 \kappa_{\infty 1}, \\
\ord_{T=0}H'_Y(T,T^m\tilY_{\infty,1})\ge \kappa_{\infty j}.
\end{gathered}
\end{align}
The final two equations are 
\begin{equation}
\label{euniV}
\Gamma_{\infty 1,-m}=1, \qquad \Gamma_{\infty 1\,0}=0.
\end{equation}

The following statement is immediate in view of the definitions and properties from Sections~\ref{spui} and~\ref{spuinf}.

\begin{proposition}
Vector~$\ph$ belongs to the set~$V$. \qed
\end{proposition}

\section{The Algebraic Set~$W$}
We write 
$$
W=W_1\cup W_2\cup W_3\cup W_4\cup W_5\cup W_6,
$$
where the sets $W_1,\ldots,W_6$ are defined below.
 
The set~$W_1$ is defined by
${\Delta=0}$.
Next, put
$$
W_2=\bigcup_{\genfrac{}{}{0pt}{}{1\le i\le \mu}{1\le j\le \nu}}W_2^{(ij)}, \qquad W_3=\bigcup_{1\le i<j\le\nu}W_3^{(ij)},
$$
where~$W_2^{(ij)}$ is defined by ${\Alpha_i=\Beta_j}$ and~$W_3^{(ij)}$ is defined by ${\Beta_i=\Beta_j}$.
  
Further, we put 
$$
W_4=\bigcup_{\genfrac{}{}{0pt}{}{i\in\{1, \ldots, \nu,\infty\}}{1\le j\le  \ell_i}}W_4^{(ij)},
$$
where the set~$W_4^{(ij)}$ is defined by the relations
\begin{align}
\label{enonvander}
\ord_{X=\Beta_i}F'_Y(X,\tilY_{ij})&> \kappa_{ij},  &&\text{when $i\ne \infty$},\\
\ord_{T=0}G'_Y(T,\tilY_{\infty j})&> \kappa_{\infty j}, &&\text{when $i=\infty$ and ${j\ne 1}$},\\
\ord_{T=0}H'_Y(T,T^m\tilY_{\infty 1})&> \kappa_{\infty j}, &&\text{when $(i,j)=(\infty,1)$}. 
\end{align}

Further, we put
$$
W_5=\left(\bigcup_{\genfrac{}{}{0pt}{}{1 \le i \le \nu}{1 \le j_1 < j_2 \le \ell_i}} W_5^{\left(i j_1 j_2\right)}\right) \cup \left(\bigcup_{2\le j_1 < j_2 \le \ell_\infty} W_5^{\left(\infty j_1 j_2\right)}\right),
$$
where $W_5^{\left(i j_1 j_2\right)}$ is defined by $\Gamma_{i j_1 \lambda\left(i j_1 j_2\right)}=\Gamma_{i j_2 \lambda\left(i j_1 j_2\right)}$ and $W_5^{\left(\infty j_1 j_2\right)}$ is defined by $\Gamma_{\infty j_1 \lambda\left(\infty j_1 j_2\right)}=\Gamma_{\infty j_2 \lambda\left(\infty j_1 j_2\right)}$, the numbers $\lambda(i,j_1,j_2)$ being defined at the end of Sections~\ref{spui} and~\ref{spuinf}.

Finally, Lemma~\ref{lber} implies that there is a proper Zariski-closed subset~$W_6$ of~$V$ such that ${\ph\notin W_6}$ and for any ${\hatphi=\left(\hatutheta, \hatualpha, \hatubeta,\hatugamma,\hatdelta\right)\in V\setminus W_6}$ the polynomial 
\begin{equation}
\label{epolhat}
Y^n + \sum_{j=0}^{n-1}\sum_{i=0}^m \widehat{\theta}_{ij}X^iY^j
\end{equation}
is irreducible and has the following property. Let~$\hatx$ and~$\haty$ be the coordinate functions on the curve~$\hatCC$ defined by~(\ref{epolhat}).  Then the effective divisor $(\haty)_\infty$ satisfies 
${\dim\LL\bigl((\haty)_\infty\bigr)=2}$.

The following statement is again immediate.

\begin{propo}
The vector~$\ph$ does not belong to the set~$W$. \qed
\end{propo}

\section{Finiteness of ${V\setminus W}$}

\label{sfin}

Here we prove that the set ${V\setminus W}$ is finite.
Let ${\hatphi=\left(\hatutheta, \hatualpha, \hatubeta,\hatugamma,\hatdelta\right)}$ be a point in ${V\setminus W}$. Then ${\hatualpha=\ualpha}$ because of~(\ref{eramV}). 

Put 
$$
\hatf(X,Y)=Y^n + \sum_{j=0}^{n-1}\sum_{i=0}^m \widehat{\theta}_{ij}X^iY^j.
$$ 
It is a $\bar\Q$-irreducible polynomial (because ${\hatphi\notin W_6}$) and defines an algebraic curve~$\hatCC$ together with rational functions ${\hatx,\haty\in \bar\Q(\hatCC)}$ satisfying ${\hatf(\hatx,\haty)=0}$. Notice that this implies that~$\haty$ is integral over $\bar\Q[\hatx]$.

Let $\hatd(X)$ be the $Y$-discriminant of $\hatf(X,Y)$. Then
$$
\hatd(X)=\hatdelta\prod_{i=1}^\mu\left(X-\alpha_i\right)^{\sigma_i}\prod_{i=1}^\nu\left(X-\hatbeta_i\right)^{\tau_i}
$$
because~$\hatphi$ satisfies~(\ref{edisV}).  Since ${\hatphi\notin W_2\cup W_3}$, the numbers~$\hatbeta_i$ are pairwise distinct and also are distinct from every~$\alpha_i$.  

The covering ${\hatCC\stackrel\hatx\to \PP^1}$  can be ramified only over the roots of $\hatd(X)$, and, perhaps, over infinity. We want to show that~$\hatx$ is unramified over the numbers~$\hatbeta_i$ and over infinity.

Fix a root~$\hatbeta_i$ and  
define
\begin{equation}
\tilhaty_{ij}(X) =\sum_{k=0}^{\kappa_{ij}}\hatgamma_{ijk}(X-\hatbeta_i)^k \qquad (j=1, \ldots, \ell_i).
\end{equation}
Then 
$$
\ord_{\hatbeta_i} \hatf(X,\tilhaty_{ij}) > 2\kappa_{ij}, \qquad \ord_{\hatbeta_i} \hatf'_Y(X,\tilhaty_{ij}) = \kappa_{ij},
$$
because~$\hatphi$ satisfies~(\ref{eserV}) and does not satisfy~(\ref{enonvander}).
Also, none of~$\tilhaty_{ij}$ is an initial segment of another, because ${\hatphi\notin W_5}$. 

Using Lemma~\ref{lhens}, we find~$\ell_i$ pairwise distinct Puiseux expansions  
$$
\haty_{i1}, \ldots, \haty_{i\ell_i}\in \bar\Q[[X-\hatbeta_i]]
$$
of~$\haty$ at~$\hatbeta_i$. satisfying  ${\ord_{\hatbeta_i} \hatf'_Y(X,\haty_{ij}) = \kappa_{ij}}$.  Since
$$
\sum_{j=1}^{\ell_i} \ord_{\hatbeta_i} \hatf'_Y(X,\haty_{ij}) = \sum_{j=1}^{\ell_i}\kappa_{ij}=\tau_i=\ord_{\hatbeta_i}\hatd(X),
$$
by~(\ref{eskij}), Lemma~\ref{lunram} implies that all~$n$ Puiseux expansions of~$\hatx$ at~$\hatbeta_i$ are in ${\bar\Q[[X-\hatbeta_i]]}$, which means that~$\hatx$ is unramified over~$\hatbeta_i$.

In a similar way we prove that~$\hatx$ is unramified over infinity (here ${\ell_\infty=n}$ and we do not need Lemma~\ref{lunram}). Moreover,~$\haty$ has at infinity ${n-1}$ Puiseux expansions without negative powers and one expansion starting from the term of degree~${-m}$. Since~$\haty$ is integral over $\bar\Q[\hatx]$, we have ${(\haty)_\infty=m\hatP}$, where~$\hatP$ is a pole of~$\hatx$. Since ${\hatphi\notin W_6}$, we have 
${\dim\LL(m\hatP)=2}$.

Thus, each ${\hatphi\in V\setminus W}$ gives rise to a pair $(\hatCC,\hatx)$, where~$\hatCC$ is an algebraic curve and~$\hatx$ an rational function on~$\hatCC$ of degree~$n$, unramified outside the points~$\alpha_i$. By Lemma~\ref{lfin}, there is only finitely many possibilities for $(\hatCC,\hatx)$. Fix one. Since  ${\dim\LL(m\hatP)=2}$,  the function~$\haty$ is uniquely defined by the equations~(\ref{euniV}). It follows that the polynomial~$\hatf$ is uniquely defined as well. Hence so is~$\hatdelta$, and the  vector~$\hatubeta$ is uniquely defined up to ordering its components. Having this order fixed, we find that~$\hatugamma$ is uniquely defined.

This proves that the set ${V\setminus W}$ is finite.

\section{Estimating the Equations Defining~$V$}

In this section we estimate the degrees and the heights of the equations defining the algebraic set~$V$. 

Since ${\kappa_{ij}\le \deg d(X)\le 2m(n-1)}$, equations defined by~(\ref{eserV}) are of degree at most 
$$
n\bigl(2m(n-1)+1\bigr)+1\le 2mn^2.
$$
Here the ``$1$'' inside the parentheses is the degree of $\tilY_{ij}$ in~$\uGamma$, and the ``$1$'' outside the parentheses is the degree of~$F$ (and of~$F'_Y$) in~$\uTheta$. 

A straightforward verification shows that the degrees of the other equations are bounded by $2mn^2$ as well. 

Now let us estimate the heights of the equations. The heights of the $\mu$ equations~(\ref{eramV}) are obviously bounded by ${h=\max\{\height(\alpha_1), \ldots, \height(\alpha_\mu)\}}$.

Estimating the heights of the remaining equations can be done with Lemma~\ref{lsombrapol}. 
All of the polynomials occurring below have rational integer coefficients. We call the  \textsl{size} of a polynomial~$p$ with coefficients in~$\Z$ (denoted by $\|p\|$) the sup-norm of the vector of its coefficients. For a non-zero polynomial~$p$ we have ${\height (p) \le \log \|p\|}$, with equality if the coefficients are co-prime. In particular, ${\height(p)=0}$ if~$p$ is of size~$1$, which is the case for many polynomials below. 

The left-hand side of~(\ref{edisV}) is a determinant 
of order ${2n-1}$ whose entries are polynomials in ${n(m+1)+1}$ variables~$X$ and~$\uTheta$, each entry being of degree at most ${m+1}$ and of size at most~$n$. Hence its height can be estimated using Lemma~\ref{lsombramat}:
$$
\height(D) \le (2n-1)\Bigl(\log n + \log (2n-1)+ (m+1)\log\bigl(n(m+1)+2\bigr)\Bigr) \le 10(mn)^2.
$$
The right-hand side of~(\ref{edisV}) is a product of at most ${2m(n-1)}$ polynomials of degree~$1$ and size~$1$ in ${\mu+\nu+1}$ variables~$\uAlpha$,~$\uBeta$ and~$X$. 
Lemma \ref{lsombrapol}~(\ref{iprod}) allows us to estimate the height of the right-hand side by the quantity 
${2m(n-1)\log(\nu+\mu+1) \le 5(mn)^2}$. 
We thereby bound the heights of the equations coming from~(\ref{edisV}) by ${10(mn)^2}$. 

Equations~(\ref{euniV}) are, obviously, of height~$0$. The height of equations coming from~(\ref{eserV}),~(\ref{eserVinf1}) and~(\ref{eserVinf2}) can be estimated using Lemma~\ref{lsombrapol}~(\ref{ig}). 
For ${i\ne \infty}$ the polynomial~$\tilY_{ij}$ is in ${\kappa_{ij}+ 2\le 2mn}$ variables~$X$,~$\Beta_j$,  $\uGamma_{ij}$.  It is of degree ${\kappa_{ij}+1\le 2mn-1}$ and of size bounded by ${2^{\kappa_{ij}}\le 4^{mn}}$. Lemma~\ref{lsombrapol}~(\ref{ig}) together with Remark~\ref{rsost} bound  the height of the polynomials $F(X,\tilY_{i,j})$ and $F'_Y(X,\tilY_{i,j})$ are bounded by the quantities 
$$
\bigl(mn\log 4+ \log 2+ 2mn\log(2mn+1)\bigr)(m+n) 
$$
and 
$$
\log n+\bigl(mn\log 4+ \log 2+ 2mn\log(2mn+1)\bigr)(m+n-1),
$$
respectively. Both do not exceed $6(mn)^3$, which bounds the heights of equations coming from~(\ref{eserV}). 
Similarly, one bounds by $12(mn)^3$ the heights of equations coming from~(\ref{eserVinf1}) and~(\ref{eserVinf2}). 

Finally, we summarize all these calculations with the following proposition.

\begin{proposition}
The algebraic set~$V$ is defined by equations of degree bounded by $2mn^2$ and height bounded by ${h+12(mn)^3}$. 
\end{proposition}


\section{The Height of~$\ph$ and the Field $\K(\ph)$}
Now we may apply Proposition~\ref{pbound}, or, more precisely, Corollary~\ref{cbound} to bound the height of the vector~$\ph$, and the number field generated by its coordinates. Recall that~$\ph$ belongs to $\bar\Q^\Omega$, where the dimension~$\Omega$ satisfies 
$$
\Omega \le 10mn+2n-7,
$$
see~(\ref{eOmega}). 
If we define~$\nabla$ and~$\Sigma$ as in Proposition~\ref{pbound}, we would have
$$
\height(f)\le \height(\ph) \le \nabla\Sigma\bigl(h+12(mn)^3\bigr)+ 2\nabla\Omega\log(\Omega+1).
$$
Furthermore, the field ${\L=\K(\ph)}$ satisfies ${[\L:\K]\le \nabla}$ and 
$$
\frac{\norm_{\L/\K}\DD_{\L/K}}{[\L:\Q]} \le 2\nabla\Sigma\bigl(h+12(mn)^3\bigr)+ 5\nabla\Omega\log(\Omega+1).
$$
Since the degrees of the equations defining~$V$ are bounded by $2mn^2$, we have   
$$
\nabla\le (2mn^2)^\Omega \le (2mn^2)^{10mn+2n-7}.
$$
Obviously, ${\Sigma\le \Omega\le 12mn}$. After trivial calculations we obtain
\begin{equation}
\label{elam'}
\height(f) \le \Lambda' (h+1), \qquad
[\L:\K]\le \Lambda', \qquad \frac{\log\norm_{\K/\Q}\DD_{\L/\K}}{[\L:\Q]}\le \Lambda' (h+1)
\end{equation}
with ${\Lambda'=(2mn^2)^{10mn+2n-3}}$. Since ${m=\genus+1}$, this proves Theorem~\ref{tmain} in the case when there is no ramified points and no Weierstrass points among the poles of~$x$.

\section{The General Case}
\label{sram}
We no longer assume that the set of poles of~$x$ has no Weierstrass and no ramified  points (called \textsl{bad} points in the sequel). Since there exists at most ${\genus^3-\genus}$ Weierstrass points and at most $2\genus$ ramified points, there exists ${\rho\in \Z}$, satisfying 
$$
|\rho|\le \genus^3+\genus\le m^3
$$
(recall that ${m=\genus +1}$) such that the fiber of~$x$ above~$\rho$ contains no bad  points. It follows that the function ${\check x=(x-\rho)^{-1}}$ has no bad points among its poles, and the previous argument applies to it. We find a number field~$\L$, a rational function ${y\in \L(\CC)}$ such that ${\L(\CC)=\L(\check x, y)}$, and a polynomial ${\check f(X,Y)\in \L[X,Y]}$ such that ${\check f(\check x, y)=0}$, 
$$
\deg_X\check f =m=\genus+1, \qquad\deg_Y\check f =n,
$$
and~(\ref{elam'}) holds with~$f$ replaced by~$\check f$ and~$h$ replaced by 
$$
\check h:=\max \Bigl\{\height \bigl((\alpha_1-\rho)^{-1}\bigr), \ldots, \height \bigl((\alpha_\mu-\rho)^{-1}\bigr)\Bigr\}. 
$$
Obviously
$$
\check h\le h+ \log \bigl(2\max\{1,|\rho|\}\bigr)\le h+ 3\log(2m),
$$
which proves~(\ref{efield}) after a short calculation. Further, the polynomial 
$$
f(X,Y):= (X-\rho)^m\check f\bigl((X-\rho)^{-1}, Y)
$$
satisfies ${f(x,y)=0}$ and 
$$
\height(f)\le \height(\check f)+3m\log(2m)
$$
by Lemma~\ref{lrho}. Again a trivial calculation  implies~(\ref{epoly}).  Theorem~\ref{tmain} is completely proved. 
\qed

\section{On the Work of Zverovich}
\label{szv}

As we already indicated in the introduction, the prototype of our proof is the work of Zverovich~\cite{Zv87}. Given a covering ${\CC\stackrel x\to\PP^1}$ and a point ${\alpha \in \PP^1}$, call the \textsl{total ramification} of~$x$ at~$\alpha$ the quantity 
$$
e(\alpha)=e_x(\alpha)=(e_1-1)+\cdots + (e_s-1),
$$ 
where ${e_1, \ldots, e_s}$ are the ramification indices of~$x$ over~$\alpha$. If particular, ${e(\alpha)>0}$ if and only if~$x$ is ramified over~$\alpha$. 

Loosely, Zverovich's argument is as follows. He defines~$x$,~$y$ and the polynomial~$f$ in (almost) the same way as we do. Then, denoting by $d(X)$ the $Y$-discriminant of~$f$, one has the equality
$$
d(X)= \prod_{i=1}^\mu(X-\alpha_i)^{e(\alpha_i)} \psi(X)^2,
$$
where~$\psi$ is a polynomial. 
Zverovich considers  the equations which follow from the relation 
\begin{equation}
\label{ezv}
D(X)= \prod_{i=1}^\mu(X-\alpha_i)^{e(\alpha_i)} \Psi(X)^2,
\end{equation}
where the unknown are the coefficients of variable polynomials~$F$ and~$\Psi$, and, as in our argument, $D(X)$ is the $Y$-discriminant of the variable polynomial~$F$. He  adds to this two equations similar to our normalization equations~(\ref{euniV}). He observes that $(f,\psi)$ satisfies his system of equations, and wants to prove that the system has finitely many equations. 

Unfortunately, Zverovich's proof of finiteness seems to be incomplete. In fact, he implicitly assumes that, for any solution $(\hat f,\widehat\psi)$ of~(\ref{ezv}), the curve $\hat\CC$, defined by ${\hat f(X,Y)=0}$,  is ramified over the points ${\alpha_1, \ldots, \alpha_\mu}$, and, moreover, the total ramification is the same as for our curve. If this were true, then Zverovich would have correctly proved that there is no other ramification, and Lemma~\ref{lfin} would imply finiteness. The problem is that a curve defined by a polynomial satisfying Zverovich's equations is not obliged \textsl{a priori} to have the same ramification at the points ${\alpha_1, \ldots, \alpha_\mu}$, as our curve, and without this  his argument does not seem to work.

We failed to repair Zverovich's argument and had to re-invent another system of equations defining our polynomial~$f$, which is much more complicated than  his one. It would be interesting to re-consider his work and try to justify his argument. This would not only improve on the estimates of this article, but would also probably imply a relatively practical  algorithm (see~\cite{DZ00} for some indications) for actual calculation of the polynomial~$f$. Evidently, our equations are too bulky for this purpose.

\subsection*{Acknowledgements}
We thank Anna Cadoret, Pierre Dèbes, Bas Edixhoven, Carlo Gasbarri and Martin Sombra for helpful discussions.


\begin{thebibliography}{99}




\normalsize



\bibitem{BSS10}
\textsc{Yu.~Bilu, M.~Strambi, A.~Surroca}, Quantitative Chevalley-Weil Theorem for Curves, in preparation. 

\bibitem{BGS94}
\textsc{J.-B.~Bost, H.~Gillet, C.~Soul\'e}, 
Heights of projective varieties and positive Green forms,
\textit{J.~Amer. Math. Soc.}~\textbf{7} (1994), 903--1027. 

\bibitem{Bi93}
\textsc{Yu. Bilu},
Effective Analysis of Integral Points on Algebraic Curves,
Ph.~D. Thesis, Beer Sheva, 2003. 

\bibitem{Deb01}
\textsc{P.~Dèbes}, 
Méthodes topologiques et analytiques en théorie inverse de Galois: théorème d'existence de Riemann, in~\cite{De01}, pp. 27--41. 

\bibitem{De01}
\textsc{B. Deschamps} (ed.), 
\textit{Arithmétique des revètements algébriques: 
Proc.  colloq. Saint-\'Etienne, March 24--26, 2000}, Séminaires et Congrès  \textbf{5},  SMF, Paris, 2001.

\bibitem{DZ00}
\textsc{O.~B.~Dolgopolova, \'E.~I.~Zverovich},  Explicit construction of global uniformization of an algebraic correspondence (Russian),  \textit{Sibirsk. Mat. Zh.}~\textbf{41}  (2000),   72--87, ii  (translated in \textit{Siberian Math.~J.}~\textbf{41}  (2000),   61--73). 

\bibitem{EJS08}
\textsc{B.~Edixhoven, R.~de Jong, J.~Schepers}, Covers of surfaces with fixed branch locus, 	\textsf{arXiv}:0807.0184. 


\bibitem{Ha93}
\textsc{R.~Hartshorne},
\textit{Algebraic Geometry}, Graduate Texts in Math.~\textbf{52}, Springer, New York, 1977. 

\bibitem{KPS01}
\textsc{T.~Krick, L.M.~Pardo, M.~Sombra},
Sharp estimates for the arithmetic Nullstellensatz,
\textit{Duke Math. J.} \textbf{109} (2001),  521--598.

\bibitem{Ph91}
\textsc{P. Philippon}, Sur des hauteurs alternatives, I, \textit{Math. Ann}.~\textbf{289} (1991),  255--283; II, \textit{Ann. Inst.
Fourier}~\textbf{44} (1994), 1043--1065; III, \textit{J. Math. Pures Appl.}~\textbf{74} (1995),  345--365.

\bibitem{Si84}
\textsc{J.~H.~Silverman}, Lower bounds for height functions, \textit{Duke Math. J}.~\textbf{51} (1984), 395--403.

\bibitem{vW93}
\textsc{B.~L.~van~der~Waerden}, \textit{Algebra II}, 6th German edition,  Springer, 1993.

\bibitem{Zv87}
\textsc{\'E.~I.~Zverovich}, 
An algebraic method for constructing the basic functionals of a Riemann surface  given in the form of a finite covering of a sphere (Russian),
\textit{Sibirsk. Mat. Zh.}~\textbf{28} (1987),  32--43, 217; (translated in \textit{Siberian Math.~J.}~\textbf{28}  (1987),   889--898). 




\end{thebibliography}
\end{document}